\documentclass[graybox]{SNmult}

\usepackage{type1cm}        % activate if the above 3 fonts are
\usepackage{makeidx}         % allows index generation
\usepackage{graphicx}        % standard LaTeX graphics tool
\usepackage{multicol}        % used for the two-column index
\usepackage[bottom]{footmisc}% places footnotes at page bottom

\usepackage{newtxtext}       % 
\usepackage[varvw]{newtxmath}       % selects Times Roman as basic font

\makeindex             % used for the subject index
\usepackage[ansinew]{inputenc}
\usepackage[T1]{fontenc}
\usepackage{float}
\usepackage{tikz}
\usetikzlibrary{arrows,matrix,fit,shapes.geometric}
\usepackage{caption}
\usepackage{ulem}
\usepackage{subcaption}
\DeclareMathOperator{\supp}{supp}
\newcommand{\bs}[1]{\boldsymbol{#1}}
\newcommand{\etal}{\textit{et al.}}

\newcommand{\figureScale}{0.917\textwidth}
\newcommand{\figureScalee}{0.75\textwidth}

\begin{document}
%%%%%%%%%%%%%
%%% Title %%%
%%%%%%%%%%%%%
\title*{Overlapping Domain Decomposition for Meshless Finite Difference Methods}
%%%%%%%%%%%%%%%
%%% Authors %%% 
%%%%%%%%%%%%%%%
\author{Alexander Westermann\orcidID{0009-0000-8452-9400}\and
Oleg Davydov\orcidID{0000-0001-8813-9485}\and 
Stefan Turek\orcidID{0000-0002-9740-6087}}
\institute{Alexander Westermann
\at Department of Mathematics, JLU Giessen, \email{alexander.westermann@math.uni-giessen.de}
\and 
Oleg Davydov
\at Department of Mathematics, JLU Giessen, \email{oleg.davydov@math.uni-giessen.de}
\and
Stefan Turek
\at Institute for Applied Mathematics, TU Dortmund University, \email{Stefan.Turek@math.tu-dortmund.de}}
\maketitle

%%%%%%%%%%%%%%%%% 
%%%% Abstract %%%
%%%%%%%%%%%%%%%%%
\abstract*{Schwarz type domain decomposition methods generally require a
partition of unity to combine solutions on subdomains.
However, in mesh-based methods it is common to organize subdomains with
minimal overlap, if any, which is facilitated by the availability
of a mesh. This study analyzes how the continuity of the partition of
unity affects the algebraic Schwarz method for Poisson and Stokes equations
from a meshless point of view, whereby the underlying differential operators are
discretized using the radial basis function finite difference (RBF-FD)
method. We demonstrate numerically that, in this setting, small overlaps improve
the performance of the domain decomposition, leading to smaller iteration counts, and
therefore no disjoint partitioning technique is required.
%%%%%%%%%%%%%%%%%
%%% Key words %%% 
%%%%%%%%%%%%%%%%%
}
\abstract{Schwarz type domain decomposition methods generally require a
partition of unity to combine solutions on subdomains.
However, in mesh-based methods it is common to organize subdomains with
minimal overlap, if any, which is facilitated by the availability
of a mesh. This study analyzes how the continuity of the partition of
unity affects the algebraic Schwarz method for Poisson and Stokes equations
from a meshless point of view, whereby the underlying differential operators are
discretized using the radial basis function finite difference (RBF-FD)
method. We demonstrate numerically that, in this setting, small overlaps improve
the performance of the domain decomposition, leading to smaller iteration counts, and
therefore no disjoint partitioning technique is required.
%%%%%%%%%%%%%%%%%
%%% Key words %%% 
%%%%%%%%%%%%%%%%%
}

%%%%%%%%%%%%%%%%%%%
%%% Introduction %%% 
%%%%%%%%%%%%%%%%%%%%
\section{Introduction}
Schwarz's domain decomposition going back to \cite{schwarz1870original} is used in modern mesh-based discretization methods 
either as a stand-alone iterative method for solving linear systems, or  as a preconditioner. In particular,  
it is applied to saddle point problems, such as the Stokes or Navier-Stokes equations, for which no classical iterative
methods can be used due to the matrix structure~\cite{gander2008schwarz,toselli2004domain}. 

Furthermore, numerous papers on finite elements, 
for example~\cite{gander2008schwarz,john2001higher,wobker2009numerical},
have discussed the application of the additive Schwarz method as a smoother in  multigrid methods, with a focus on
decomposing into very small subdomains (batch- or cell-based). These methods originate from the work of
Vanka~\cite{vanka1986block}, which employed a symmetric coupled Gauss-Seidel smoother 
that can be identified as an algebraic Schwarz method. 

In this paper we explore the application of the additive Schwarz method to linear systems arising from a meshless finite
difference method on irregular nodes. The main question is how to design subdomains without the help of a mesh, and we
demonstrate that a crucial parameter is the amount of the overlap between different subdomains. In mesh-based methods it is
common~\cite{gander2018does} to use subdomains that build a disjoint partition of the full computation domain. However,
our experiments show that this is undesirable for the meshless methods.

The second question we address is the influence of the choice of the partition of unity (PoU) used to combine subsolutions when
subdomains overlap. It turns out that both versions of the algebraic additive Schwarz method introduced 
in~\cite{cai1999restricted},  the restricted additive Schwarz method (RAS) and the 
additive Schwarz method with harmonic extension (ASH) perform well in the meshless setting if PoU is continuous, but 
ASH typically outperforms RAS with the simpler discontinuous PoU.

In our numerical experiments we discretize a Poisson problem in 1D and 2D and a Stokes problem in 2D by using the radial
basis  function finite difference (RBF-FD) method with polyharmonic RBF and a polynomial extension, see for 
example~\cite{fornberg2015primer}. 

The paper is organized as follows. After a brief Section~\ref{Schwarz} that introduces the 
two versions of the iterative Schwarz methods RAS and ASH studied in this paper in the meshless context, 
we present in Section~\ref{tests} our numerical investigation of their performance in conjunction with two types of PoU for
three examples of boundary value problems. Section~\ref{Conclusion} provides a conclusion.

%%%%%%%%%%%%%%%%%%%%%%%
%%%% Schwarz method %%% 
%%%%%%%%%%%%%%%%%%%%%%%
\section{Schwarz methods}\label{Schwarz}

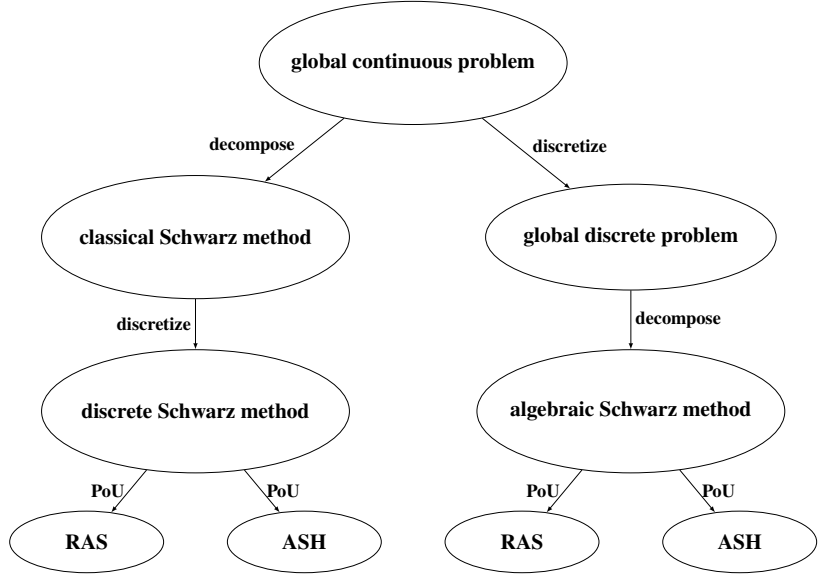
\begin{figure}[H]
	\centering
	\resizebox{0.88\textwidth}{!}{%
		\begin{tikzpicture}
			\begin{scope}
				\node[ellipse, draw, fit={(0,0) (5,2)}, xshift=5cm, inner sep=0pt, 
				label=center:\textbf{\Large{global continuous problem}}, thin] (A) {};
				\node[ellipse, draw, fit={(0,0) (5,2)}, xshift=10cm, yshift=-4cm, inner sep=-4pt, 
				label=center:\textbf{\Large{global discrete problem}}, thin] (B) {};
				\node[ellipse, draw, fit={(0,0) (5,2)}, yshift=-4cm, inner sep=0pt, 
				label=center:\textbf{\Large{classical Schwarz method}}, thin] (C) {};
				\node[ellipse, draw, fit={(0,0) (5,2)}, xshift=10cm, yshift=-8cm, inner sep=0pt, 
				label=center:\textbf{\Large{algebraic Schwarz method}}, thin] (D) {};
				\node[ellipse, draw, fit={(0,0) (5,2)}, yshift=-8cm, inner sep=0pt, 
				label=center:\textbf{\Large{discrete Schwarz method}}, thin] (E) {};
				\draw[-latex] (A)--(B);
				\draw[-latex] (A)--(C);
				\draw[-latex] (B)--(D) node [pos=0.5,right] {\textbf{\large{decompose}}};
				\draw[-latex] (C)--(E) node [pos=0.5,left] {\textbf{\large{discretize}}};
				\node[ellipse, draw, fit={(0,0) (2.5,1)}, xshift=+3.75cm,  yshift=-10.5cm, inner sep=0pt, 
				label=center:\textbf{\Large{ASH}}, thin] (F) {};
				\node[ellipse, draw, fit={(0,0) (2.5,1)}, xshift=-1.25cm,  yshift=-10.5cm, inner sep=0pt, 
				label=center:\textbf{\Large{RAS}}, thin] (G) {};
				\node[ellipse, draw, fit={(0,0) (2.5,1)}, xshift=+13.75cm,  yshift=-10.5cm, inner sep=0pt, 
				label=center:\textbf{\Large{ASH}}, thin] (H) {};
				\node[ellipse, draw, fit={(0,0) (2.5,1)}, xshift=8.75cm,  yshift=-10.5cm, inner sep=0pt, 
				label=center:\textbf{\Large{RAS}}, thin] (I) {};
				\draw[-latex] (E)--(F) node [pos=0.5,right] {\textbf{\large{PoU}}};
				\draw[-latex] (E)--(G) node [pos=0.5,left] {\textbf{\large{PoU}}};
				\draw[-latex] (D)--(H) node [pos=0.5,right] {\textbf{\large{PoU}}};
				\draw[-latex] (D)--(I) node [pos=0.5,left] {\textbf{\large{PoU}}};
				\node at (11.1,-.9) {\textbf{\large{discretize}}};
				\node at (3.8,-.9) {\textbf{\large{decompose}}};
			\end{scope}
	\end{tikzpicture}}
	\caption{Family tree of Schwarz methods}
	\label{fig:familiy_tree}
\end{figure}

Various types of iterative Schwarz methods are  summarized in Fig.~\ref{fig:familiy_tree}. It starts with a
\textit{global continuous problem} on the top, for example the Poisson equation with Dirichlet boundary conditions
\begin{align*}
	\left.\begin{cases}
		-\Delta u=f & \text{in }\Omega\\
		\phantom{-\Delta}u = g & \text{on } \Gamma=\partial\Omega
	\end{cases}\right\}.
\end{align*}
The two main branches in  Fig.~\ref{fig:familiy_tree} are 
distinguished by whether decomposition or discretization is applied first.

The general idea behind every Schwarz method is as follows: If the problem is too big or the geometry is too  complex, then
several subproblems corresponding to different subdomains of $\Omega$ are solved, and these solutions are combined into an
approximation of the global problem. By iterating this process the convergence of  the global approximation is achieved. For
the sake of simplicity of the exposition we only formulate the subdivision into two subproblems on subdomains $\Omega_1$
and $\Omega_2$, with $\bar\Omega=\bar\Omega_1\cup\bar\Omega_2$.
Furthermore, only the additive  version of the Schwarz method without any optimization, such as the optimized RAS
in~\cite{St-Cyr2007optimized}, is  considered. 

By applying the decomposition to the global problem, the decomposed 
continuous problem can be formulated with a starting solution $u^0$ and $n\ge0$
\begin{align*}
	\left.\begin{cases}
		-\Delta v^{n+1}=f^{(1)} & \text{in }\Omega_1\\
		\phantom{-\Delta}v^{n+1} = g^{(1)} & \text{on } \Gamma_1=\partial\Omega_1
	\end{cases}\right\}\quad\wedge\quad
	\left.\begin{cases}
		-\Delta w^{n+1}=f^{(2)} & \text{in }\Omega_2\\
		\phantom{-\Delta}w^{n+1} = g^{(2)} & \text{on } \Gamma_2=\partial\Omega_2
	\end{cases}\right\}
\end{align*}
where $v^{n+1}$ and $w^{n+1}$ are the solutions of each subproblem, whereas $f^{(i)}=f|_{\Omega_i}$,
$g^{(i)}=u_n|_{\Gamma_i}$, $i=1,2$. 
The new global solution $u^{n+1}$ is computed using prolongation operators $P_1$ and $P_2$ that extend
$v^{n+1}$ and $w^{n+1}$ to $\Omega$, and a combining operator $C$ resulting in $u^{n+1} = C\left(P_1 v^{n+1}, P_2 w^{n+1}\right).$
This iterative algorithm is called the \textit{classical Schwarz method} in this work.
After some numerical methods are applied to solve the continuous subproblems and operators $P_i$ and $C$
are also implemented numerically, a \textit{discrete Schwarz method} 
is defined that produces a discrete approximation $u^{n+1}_h$ of $u^{n+1}$.

The second branch of the family tree starts with the discretization instead of the decomposition so that a 
\textit{global discrete problem}
\begin{align*}
	\left.\begin{cases}
		A_h u_h = f_h & \text{in }\Omega_h\\
		\phantom{A_h}u_h = g_h & \text{on } \Gamma_h
	\end{cases}\right\}
\end{align*}
is set up.  Here, the subscript $h$ symbolizes the discrete 
version of each function as a vector, $A_h$ stands for a matrix that represents 
the discrete Poisson operator, and $\Omega_h$ and $\Gamma_h$ are the discretized domain and its boundary,
given by a grid or mesh, or just a node set in the case of meshless methods. 
Similar to the classical Schwarz method, the decomposition process can be 
applied to the discrete problem by selecting matrix and vector entries corresponding to the indices of each discrete 
subdomain $\Omega_{h,i}\subset\Omega_h$.
In this paper we consider only the trivial prolongation operators,  whereby continuous solutions are extended by 
zero outside of $\Omega_i$ and the discrete vectors are filled with zeros to match the size of the global vectors.

One possible realization of combining multiple local solutions into a global one is to use a partition of unity.
For that we choose PoU functions $\varphi_i$ that satisfy  $\Omega_i\subset\supp\varphi_i$ and sum up to one globally, 
that is, in the case of two subdomains $\varphi_1+\varphi_2 \equiv 1$ in $\Omega$, leading   to the combining
operator $C\big(v, w\big)={\varphi_1}v + {\varphi_2}w$. Thus, the global solution is obtained in the continuous case as
$$
u^{n+1}={\varphi_1}P_1 v^{n+1} + {\varphi_2}P_2w^{n+1}=u^{n}+{\varphi_1}P_1 v_{\delta}^{n+1} + {\varphi_2}P_2w_{\delta}^{n+1},$$
where the second version uses the residuals $v_{\delta}^{n+1}=v^{n+1}-u^{n}$ and $w_{\delta}^{n+1}=w^{n+1}-u^{n}$.

The discrete algebraic Schwarz method is called restricted additive Schwarz (RAS)
preconditioner in \cite{cai1999restricted}. In the same paper Cai and Sarkis also introduced 
the additive Schwarz method with harmonic extension (ASH) as a second way to include the PoU 
by applying it to the right hand sides of the subproblems before solving the linear systems.
We formulate both RAS and ASH in the residual form:
 with a global initial solution vector $u_h^0$ and $n\ge0$,
\begin{itemize}
	\item \textbf{RAS:} apply the {PoU} \textbf{after} solving each subproblem.\\
	 Solve subsystems:\\
	 \resizebox{\linewidth}{!}{
 	\begin{minipage}{\linewidth}
		\begin{align*}
			\left.\begin{cases}
				A_h^{(1)}v_{\delta,h}^{n+1} = f_h^{(1)} + A_h^{(1)} u_h^n|_{\Omega_{h,1}} & \text{in }\Omega_{h,1}\\
				\phantom{A_h^{(1)}}v_{\delta,h}^{n+1} = 0 & \text{on } \Gamma_{h,1}
			\end{cases}\right\}\quad\wedge\quad
			\left.\begin{cases}
				A_h^{(2)} w_{\delta,h}^{n+1} = f_h^{(2)} + A_h^{(2)} u_h^n|_{\Omega_{h,2}} & \text{in }\Omega_{h,2}\\
				\phantom{A_h^{(2)}}w_{\delta,h}^{n+1} = 0 & \text{on } \Gamma_{h,2}
			\end{cases}\right\}.
		\end{align*}
	\end{minipage}}\\[1em]
	Obtain a global solution from the subsolutions $v_{\delta, h}^{n+1},w_{\delta, h}^{n+1}$:
	\[
		u_h^{n+1} = u_h^n + {\varphi_{h,1}}P_{h,1} v_{\delta, h}^{n+1} + {\varphi_{h,2}}P_{h,2}
		w_{\delta, h}^{n+1}.
	\]
	
	\item \textbf{ASH:} Apply the {PoU} \textbf{before} solving each subproblem:\\
	Solve subsystems:\\
	\resizebox{\linewidth}{!}{
	\begin{minipage}{\linewidth}
		\begin{align*}
			\left.\begin{cases}
				A_h^{(1)} v_{\delta,h}^{n+1} = {\varphi_{h,1}}\left(f_h^{(1)} + A_h^{(1)} 
				u_h^n|_{\Omega_{h,1}}\right) & \text{in }\Omega_{h,1}\\
				\phantom{A_h^{(1)}}v_{\delta,h}^{n+1} = 0 & \text{on } \Gamma_{h,1}
			\end{cases}\right\}\quad\wedge\quad
			\left.\begin{cases}
				A_h^{(2)}w_{\delta,h}^{n+1} = {\varphi_{h,2}}\left(f_h^{(2)} + A_h^{(2)} 
				u_h^n|_{\Omega_{h,2}}\right) & \text{in }\Omega_{h,2}\\
				\phantom{A_h^{(2)}}w_{\delta,h}^{n+1} = 0 & \text{on } \Gamma_{h,2}
			\end{cases}\right\}.
		\end{align*}
	\end{minipage}}\\[1em]
	Obtain a global solution from the subsolutions $v_{\delta, h}^{n+1},w_{\delta, h}^{n+1}$:
	\[
		u_h^{n+1} = u_h^n + P_{h,1} v_{\delta,h}^{n+1} + P_{h,2} w_{\delta,h}^{n+1}.
	\]
\end{itemize}
The iterations are repeated until the root mean square residual $P_{h,1} v_{\delta,h}^{n+1} + P_{h,2} w_{\delta,h}^{n+1}$
reduces below a given tolerance $\varepsilon>0$. We used $\varepsilon=10^{-7}$ in all numerical experiments.

%%%%%%%%%%%%%%%%%%%%%%%%%
%%% Numerical results %%%
%%%%%%%%%%%%%%%%%%%%%%%%%
\section{Numerical results}\label{tests}
In the numerical experiments we will use two types of PoU:
\begin{itemize}
	\item \textbf{dPoU:} A discontinuous PoU depending on how many subdomain include a point $x$:
	\[
	\varphi_i(x) =
	\begin{cases}
		\big(\#\{k=1,\ldots,N_{DEC}:x\in\Omega_k\}\big)^{-1}, &  x\in\Omega_i,\\
		0, & \text{otherwise},
	\end{cases}
	\]
	where $N_{DEC}$ is the number of subdomains.
	The property $\sum_i \varphi_i = 1$ is directly fulfilled.
	
	\item \textbf{cPoU:}  A continuous PoU as a bump function depending on the the choice of centers $z_i\in\Omega_i$ and
	the maximum distance 
	$r_i := \sup\{\|x-z_i\|_2:x\in\Omega_i\}$ of the subdomain $\Omega_i$:
	\[
	\tilde\varphi_i(x) =\begin{cases}
		\exp\Big(-\frac{r_i^2}{r_i^2 - \|x-z_i\|_2^2}\Big), & x\in\Omega_i,\\
		0, & \text{otherwise},
	\end{cases}
	\]
	with $\varphi_i=\tilde\varphi_i/\sum_k \tilde\varphi_k$.
\end{itemize}

Each differential operator is discretized using the RBF-FD method with the PHS kernel of the form $r^3$ and a polynomial extension of
order 3 for the first derivatives and order 4 for the second derivatives, ensuring the convergence order of approximately 2, as
discussed in~\cite{westermann2025stability} for the Stokes case. 
The sets of influence consist of the same number of nearest neighbors of a given node, in particular 
10 for the Poisson equation in 1D, and 25 for the 2D Poisson case. In the case of the Stokes equations, the numbers are
20 for the gradient and divergence, and 25 for the Laplacian.

Discretization nodes are obtained from the quasi-random Halton sequence.

All computations were performed on a machine with an Intel Xeon E5-2640 v3 CPU (2.60 GHz)
using MATLAB R2024b and basic MATLAB functions of mFDlab~\cite{mFDlab} for RBF-FD.

We use direct solvers to obtain the subsolutions because they are relatively small and we also want to exploit the important 
advantage of the algebraic Schwarz method that the submatrices remain unchanged throughout the iteration process. 
Therefore, we compute LU decompositions of the matrices $A_h^{(i)}$ only once before the iteration starts by MATLAB's
\texttt{lu} command, and solve the linear systems within the iteration loop by forward and backward substitution using
MATLAB's backslash operator.

\vspace{-10pt}

\subsection{Poisson 1D}

\normalsize
In the first test, we consider a one-dimensional Poisson problem, similar to the example provided in~\cite{efstathiou2003restricted},
with the solution $u(x) = 0.5\left(x-x^2\right)$ on the unit interval $[0,1]$. The domain is decomposed into two 
overlapping subdomains $\Omega = \Omega_1\cup\Omega_2$ where $\Omega_1 = (0,\frac23)$ and $\Omega_2 = (\frac13,1)$.

\begin{figure}[t!]
	\centering
	\begin{subfigure}{0.49\textwidth}
		\centering
		\includegraphics[width = 
		0.7\textwidth]{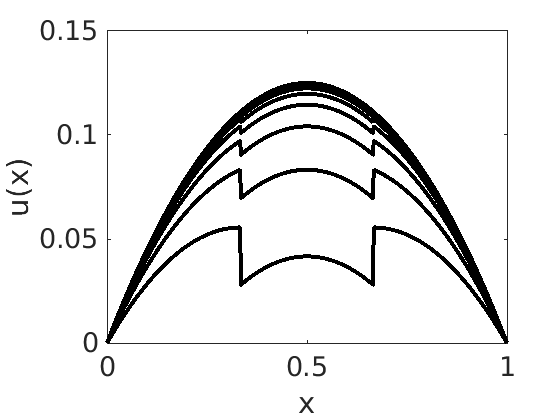}
		\caption{RAS-dPoU}
	\end{subfigure}
	\hfill
	\begin{subfigure}{0.49\textwidth}
		\centering
		\includegraphics[width = 
		0.7\textwidth]{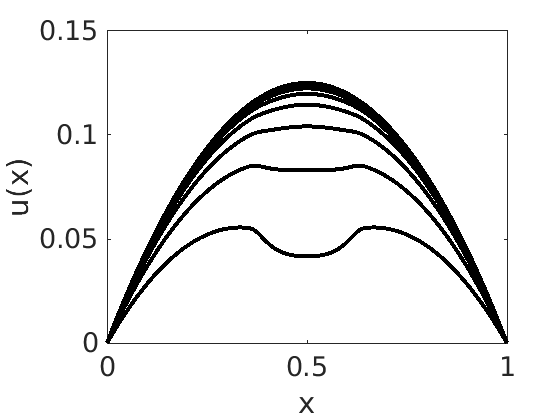}
		\caption{RAS-cPoU}
	\end{subfigure}
	\begin{subfigure}{0.49\textwidth}
		\centering
		\includegraphics[width = 
		0.7\textwidth]{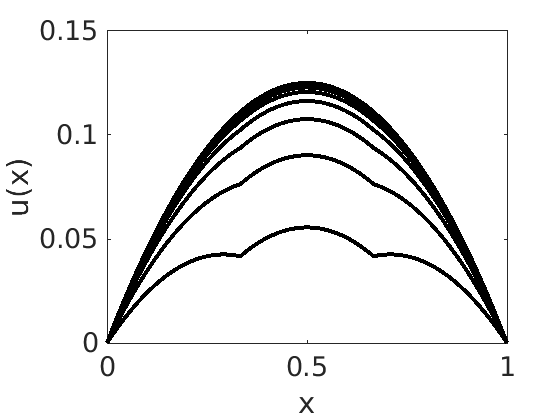}
		\caption{ASH-dPoU}
	\end{subfigure}
	\hfill
	\begin{subfigure}{0.49\textwidth}
		\centering
		\includegraphics[width = 
		0.7\textwidth]{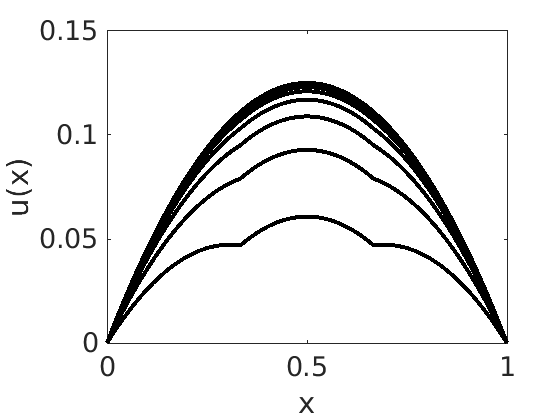}
		\caption{ASH-cPoU}
	\end{subfigure}
	\caption{Iterative solution process of algebraic Schwarz method in different combinations of RAS or ASH versions and  
		discontinuous or continuous PoU}
	\label{fig:Poisson1D_RASvsASH}
\end{figure}

Fig.~\ref{fig:Poisson1D_RASvsASH} (a)--(b) visualizes the iterative process of the global solution 
for the RAS method. It can be seen that the two different PoUs result in subsolutions of
different character. On the one hand, if a 
discontinuous PoU is applied, a discontinuous global solution is obtained. This is natural
since a discontinuous function is applied to globalizing continuous subsolutions. 
On the other hand, the application of the continuous PoU delivers a continuous 
global solution. This phenomenon disappears for ASH, where the PoU is applied before solving the subproblems. 
Moreover, as can be seen Fig.~\ref{fig:Poisson1D_RASvsASH} (c)--(d), there is no significant difference between the solutions
obtained with dPoU and cPoU.

Discontinuous subsolutions in the RAS-dPoU version have a significant drawback compared with the 
ASH-style version. First, a discontinuous approximation is clearly undesirable for a smooth solution 
of the continuous problem.
Moreover, the residuum of the discrete Poisson operator is much greater in this case because the 
discrete Poisson operator is applied to a discontinuous solution. This leads to a higher tolerance for the termination
condition, and, in  combination with a high condition number of the 
global matrix, to RAS-dPoU not converging in terms of the global residual error for finer discretizations.
We will demonstrate that even for smaller systems with lower condition numbers, RAS-dPoU 
requires more iterations than ASH-dPoU to reach a given tolerance.

To compare the performance of the methods quantitatively, and investigate the dependence on the amount of the overlap, we set
up the following experiments, both in 1D and 2D. The discretized domain $\Omega_h$ with $N_{max}$ nodes is 
decomposed into $N_{DEC}$ nearly disjoint subdomains $\Omega_{h,i}$ of approximately the same size $N_{dis}$ 
with the help of MATLAB's command \texttt{kmeans} that also returns the centroids of the subdomains. 
Each subdomain is then extended to achieve a desired number $N\le N_{max}$ of nodes, for 
$N\in[N_{min},2N_{min},3N_{min},\dots]$ by including additional nodes within an increasing distance from the centroid
until a given number $N$ of nodes is reached. Here, $N_{min}\ge N_{dis}$ is the minimum number 
of nodes required for the Schwarz method to deliver stable solutions, determined experimentally.
We refer to Figures~\ref{fig:Example_Nmin} and \ref{fig:Example_2Nmin_Stokes} for the illustrations in the 2D case.
With this design of the experiments we ensure that the number of subdomains is fixed, which  reduces the influence of
parallel computation when we compare the influence of the more or less significant overlap.
To measure the amount of the overlap  we compute the \textit{overlap ratio}
\[
	\rho(N) = \frac{N-N_{dis}}{N_{max}-N_{dis}},
\]
where $N_{dis}\le N\le N_{max}$ so that $\rho\in[0,1]$. Thus, $\rho=0$ 
stands for the disjoint decomposition and $\rho=1$ for $\Omega_i=\Omega$ for all $i\in\{1,\dots,N_{DEC}\}$. 
This overlap ratio must be interpreted as a global parameter, and it can only provide information about the relative
(local) overlap to a certain extent.

In Fig.~\ref{fig:Poisson1D_Test} and related figures below, we plot the number of
iterations $N_{iter}$ needed for the algebraic Schwarz method in RAS- and ASH-style with both discontinuous and continuous PoU to converge. In addition, in order to demonstrate how the overlap	ratio influences the overall computation cost,
we plot the total computation time in seconds for the RAS method. Clearly, this timing depends on software and hardware.
In our tests the for-loops were implemented using MATLAB's \texttt{parfor} command, with 16 parallel workers on our machine. 

We observe from the results in Fig.~\ref{fig:Poisson1D_Test}
that the overall best configurations in 1D appear to be either RAS with continuous PoU or ASH 
with discontinuous PoU. In both cases the highest possible overlap is optimal, presumably due to the small bandwidth
of the system matrix in 1D.

A relatively weak performance of RAS with discontinuous PoU
is in a good agreement with the previous observations that discontinuous solutions produce higher residual errors, 
especially at the beginning of the iteration process. 

It is surprising that the maximum overlap is optimal, even though
in this case there is no decomposition as all subdomains coincide with $\Omega$, and we in fact just repeat the same
computation $N_{DEC}$ times. However, as we will see, this does not anymore happen in 2D. A plausible explanation of this phenomenon is
that the RBF-FD system matrices for the 1D Poisson equation are sparse and banded with a low bandwidth. Therefore, the direct
method we use to solve the linear systems does not suffer from fill-in and requires just $\mathcal{O}(N)$ operations, whereas 
no iterations are needed when the global system is solved, with its numerically zero residual. 

\vspace{-10pt}

\subsection{Poisson 2D}\label{Poisson2D}

\begin{figure}[p]
\hspace{24pt}
	\begin{minipage}{0.49\textwidth}
		\includegraphics[width = \figureScalee]{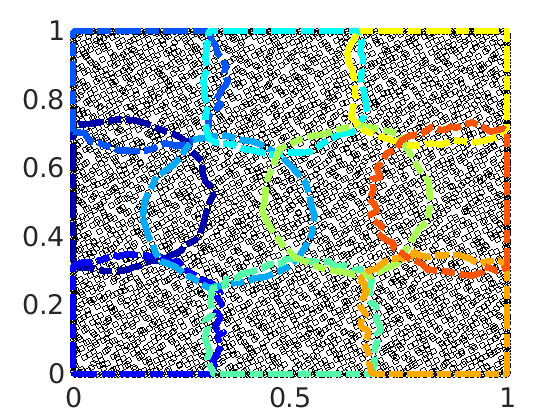}
	\end{minipage}
	\hfill
	\begin{minipage}{0.49\textwidth}
		\includegraphics[width = \figureScalee]{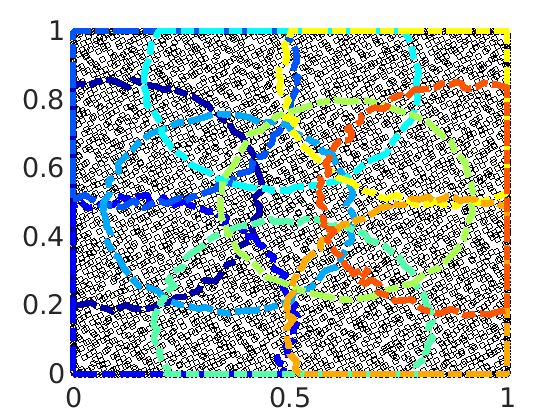}
	\end{minipage}
	\caption{Decomposition example for the Poisson problem in 2D, with $N_{DEC} = 10$ and $N_{min} = 504$ (left) and 
		$2N_{min}=1008$ (right) nodes per subdomain.}
	\label{fig:Example_Nmin}
\end{figure}

\begin{figure}[p]
	\centering
	\begin{subfigure}{0.49\textwidth}
		\centering
		\includegraphics[width = 
		\figureScale]{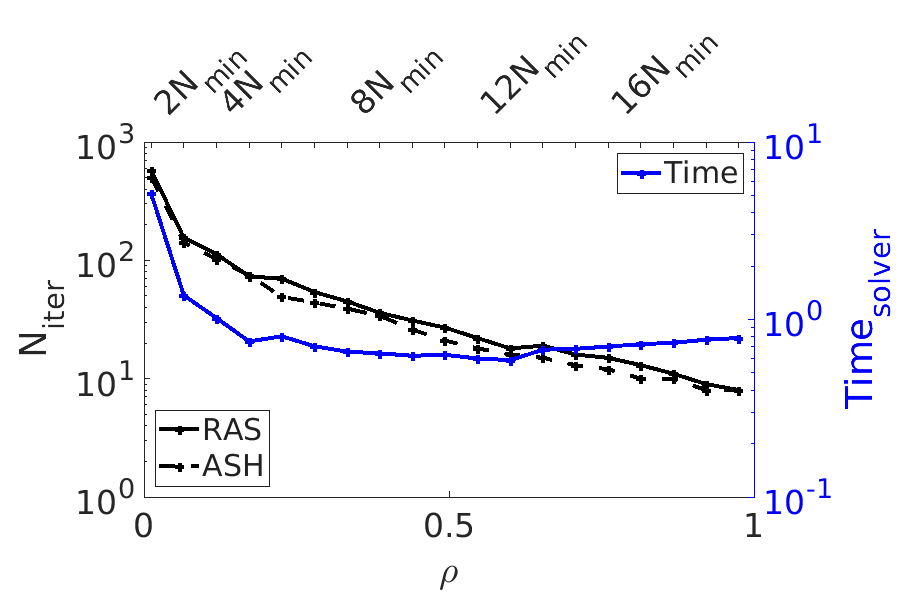}
		\caption{discontinuous PoU}
	\end{subfigure}
	\hfill
	\begin{subfigure}{0.49\textwidth}
		\centering
		\includegraphics[width = 
		\figureScale]{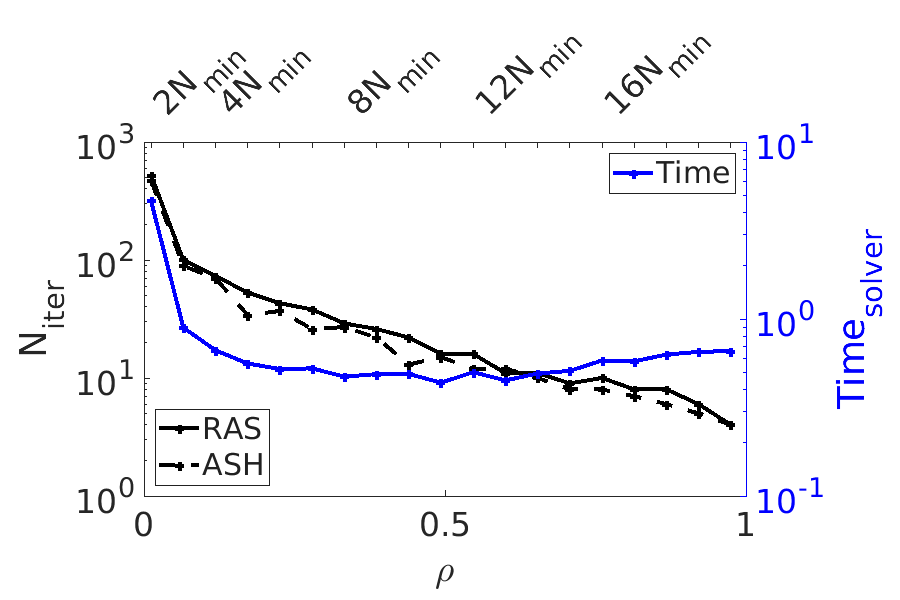}
		\caption{continuous PoU}
	\end{subfigure}
	\begin{subfigure}{0.49\textwidth}
		\centering
		\includegraphics[width = 
		\figureScale]{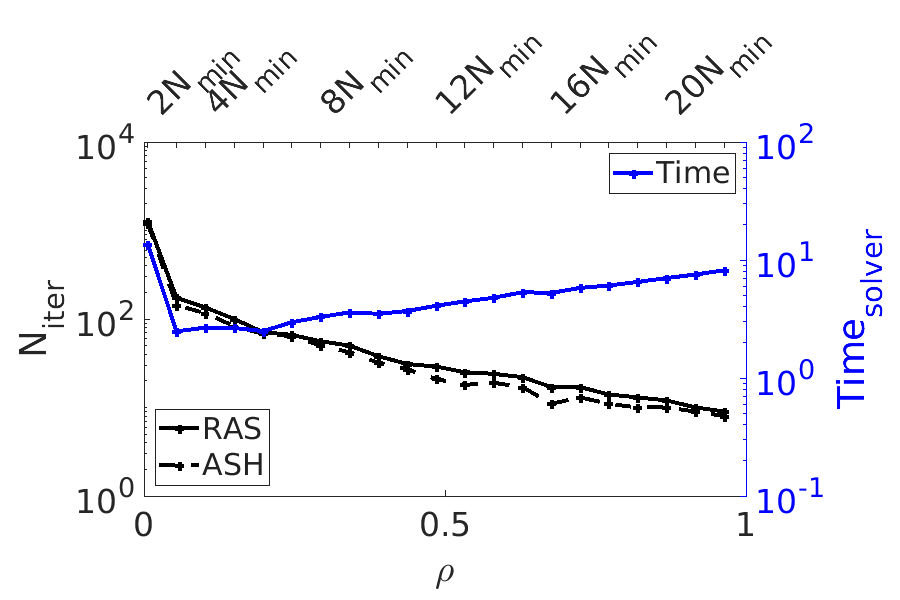}
		\caption{discontinuous PoU}
	\end{subfigure}
	\hfill
	\begin{subfigure}{0.49\textwidth}
		\centering
		\includegraphics[width = 
		\figureScale]{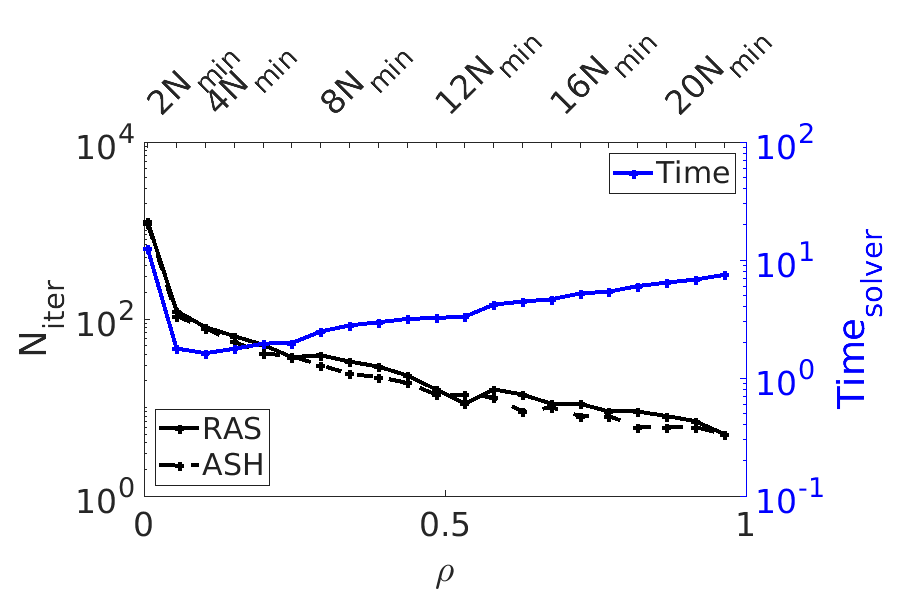}
		\caption{continuous PoU}
	\end{subfigure}
	\begin{subfigure}{0.49\textwidth}
		\centering
		\includegraphics[width = 
		\figureScale]{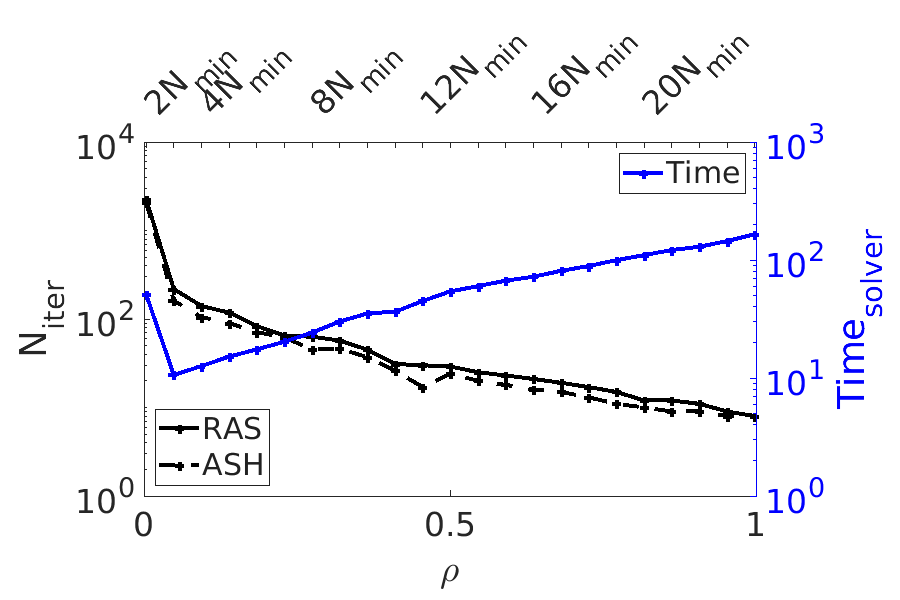}
		\caption{discontinuous PoU}
	\end{subfigure}
	\hfill
	\begin{subfigure}{0.49\textwidth}
		\centering
		\includegraphics[width = 
		\figureScale]{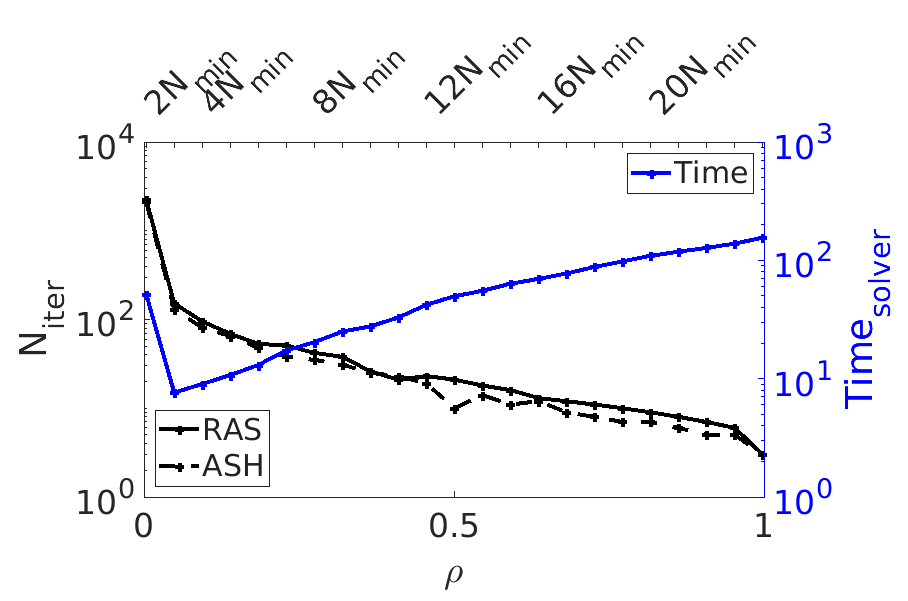}
		\caption{continuous PoU}
	\end{subfigure}
	\caption{Poisson 2D: Comparison of the number of iterations.
	(a)--(b): $N_{max}=4225$, $N_{DEC} = 25$, $N_{min}=217$, $N_{dis}=169$; 
	(c)--(d): $N_{max}=16641$, $N_{DEC} = 25$, $N_{min}=765$, $N_{dis}=666$; 
	(e)--(f): $N_{max}=66049$, $N_{DEC} = 25$, $N_{min}=2866$, $N_{dis}=2642$.
   }
	\label{fig:Poisson2D_Test}
\end{figure}

\normalsize
The test problem for this section is the Poisson equation in $[0,1]^2$ with Dirichlet boundary conditions and 
solution $u(x,y) = \sin(\pi x)\cos(\pi y)$. The setup of the experiments is the same as before. However, in contrast to 1D,
where we used $N_{min} = N_{dis}$, now  $N_{min}$ is greater 
than $N_{dis}$ because otherwise the iterations do not converge, such that some minimum overlap 
seems to be required for mere stability. We illustrate in Fig.~\ref{fig:Example_Nmin} typical
choice of subdomains for in this case. We see that the overlap is almost negligible for 
$N=N_{min}$, but already significant for  $N=2N_{min}$. Comparison of the number of iterations for RAS and ASH with dPoU and
cPoU, as well as timing is provided in Fig.~\ref{fig:Poisson2D_Test}

\begin{figure}[H]
	\begin{subfigure}{0.49\textwidth}
		\centering
		\includegraphics[width = 
		\figureScale]{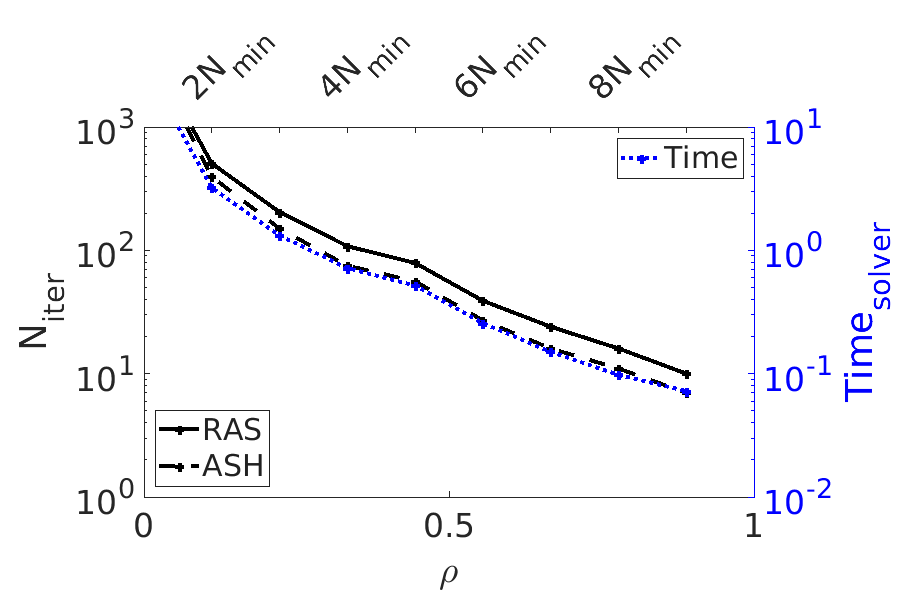}
		\caption{discontinuous PoU}
	\end{subfigure}
	\hfill
	\begin{subfigure}{0.49\textwidth}
		\centering
		\includegraphics[width = 
		\figureScale]{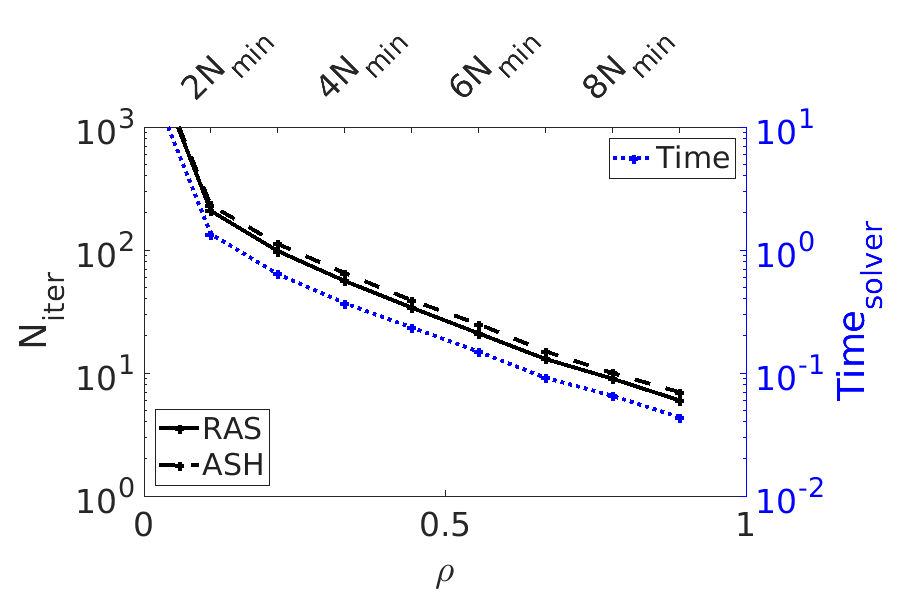}
		\caption{continuous PoU}
	\end{subfigure}
	\begin{subfigure}{0.49\textwidth}
		\centering
		\includegraphics[width = 
		\figureScale]{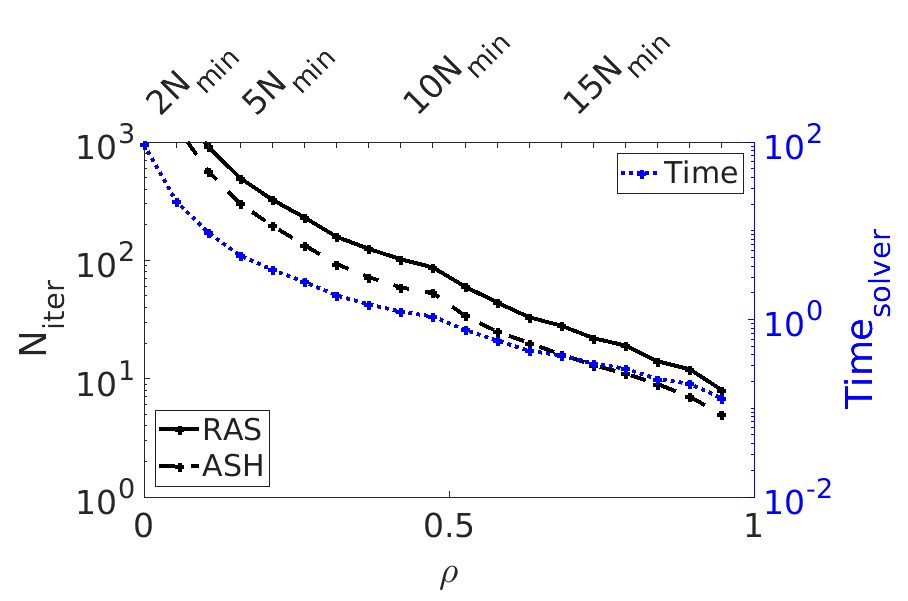}
		\caption{discontinuous PoU}
	\end{subfigure}
	\hfill
	\begin{subfigure}{0.49\textwidth}
		\centering
		\includegraphics[width = 
		\figureScale]{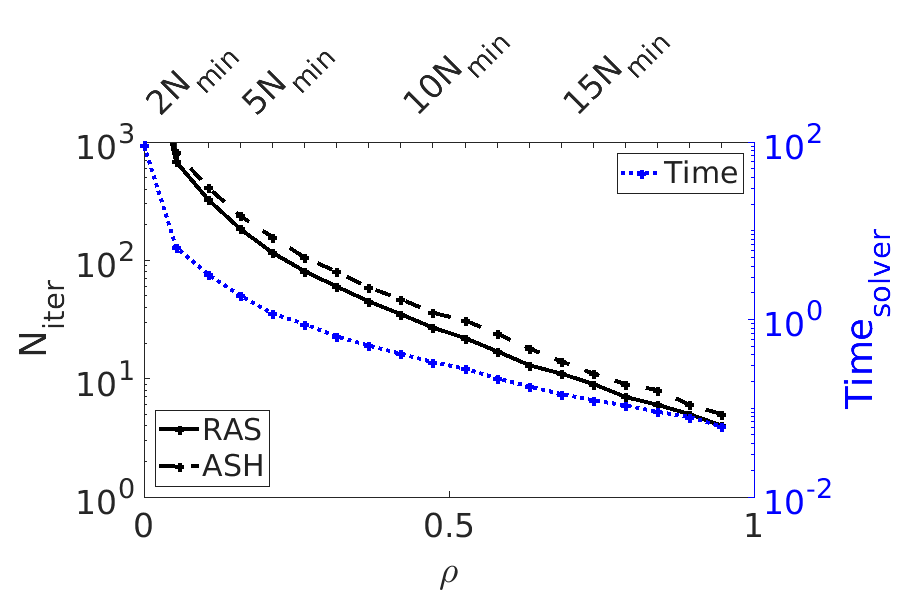}
		\caption{continuous PoU}
	\end{subfigure}
	\begin{subfigure}{0.49\textwidth}
		\centering
		\includegraphics[width = 
		\figureScale]{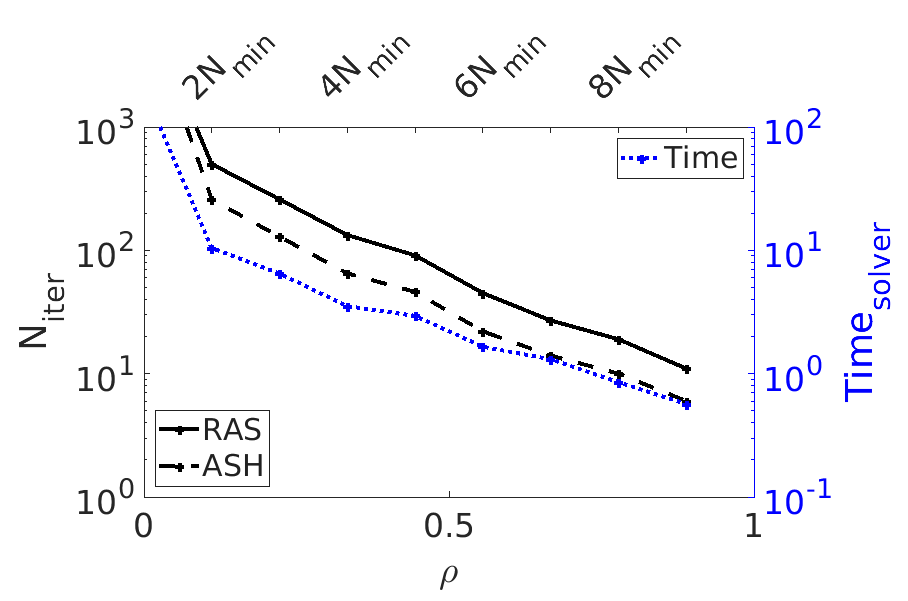}
		\caption{discontinuous PoU}
	\end{subfigure}
	\hfill
	\begin{subfigure}{0.49\textwidth}
		\centering
		\includegraphics[width = 
		\figureScale]{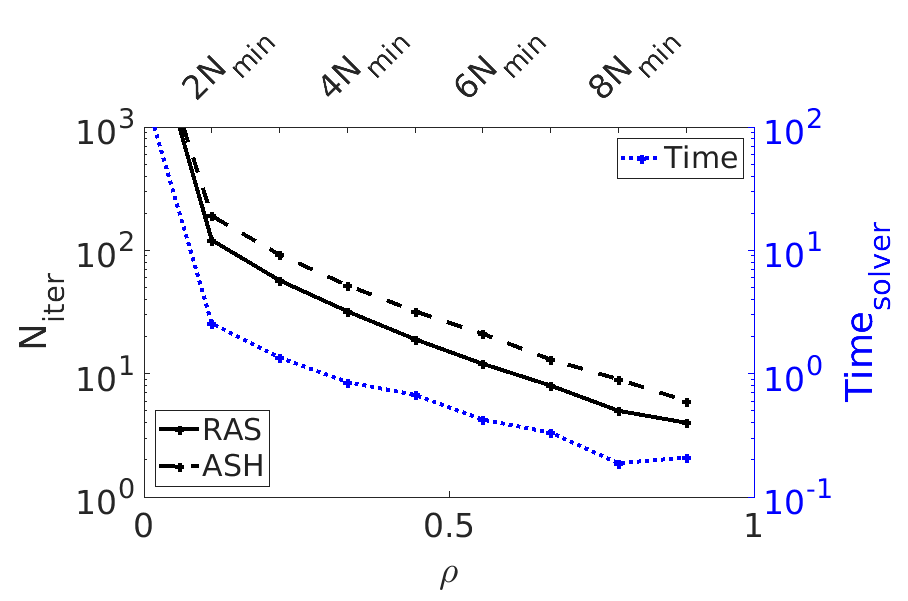}
		\caption{continuous PoU}
	\end{subfigure}
	\caption{Poisson 1D: Comparison of the number of iterations to convergence.
		(a)--(b): $N_{max}=10^3$, $N_{DEC} = 10$, $N_{min}=100$; (c)--(d): $N_{max}=10^4$, $N_{DEC} = 20$, $N_{min}=500$; 
		(e)--(f): $N_{max}=10^5$, $N_{DEC} = 10$, $N_{min}=10^4$.}
	\label{fig:Poisson1D_Test}
\end{figure}

In contrast to 1D, ASH now achieves smaller iteration numbers than RAS in almost all cases, 
even for the continuous PoU. For example, in the settings of Fig.~\ref{fig:Poisson2D_Test}(e) and (f) for $N=2N_{min}$, we have $N_{iter}(RAS)=218$ versus $N_{iter}(ASH)=161$, and $N_{iter}(RAS)=152$ 
versus $N_{iter}(ASH)=130$, respectively.
Moreover, there is a striking change in how the computation time depends on the amount of overlap.
After an initial step drop between until $N\approx 2N_{min}$ or $\rho\approx 0.05$, 
the time starts stabilizing for larger $N$, and then grows, which is more pronounced for larger sizes of the global linear system.
Clearly, this can be attributed to the higher bandwidth of the submatrices in two dimensions, and fill-in appearing in the
cause of LU-decomposition of the submatrices.

In order to investigate the cost of the method more closely, we have estimated the numerical rate of the 
growth of the computation time as function of the size $N$ of the submatrices. It turns out that the cost of computing the 
LU-decomposition of a single matrix $A_h^{(i)}$  grows approximately at the rate $\mathcal{O}(N^{1.8})$, while 
the cost of solving the system at each iteration is about  $\mathcal{O}(N^{1.4})$. This makes a big difference to the 1D case
where both rates are  $\mathcal{O}(N)$ due to the absence of the fill-in, which explains the observed behavior.

\vspace{-10pt}

\subsection{Stokes 2D}\label{Stokes2D}

\normalsize
Finally, we consider the incompressible Stokes equations in $\bar\Omega=[0,1]^2$ in the velocity-pressure formulation, $\Delta \bs u + \nabla p = \bs f$ with $\nabla\cdot \bs u = 0$ in $\Omega$, and boundary conditions $\bs u = \bs g$ on $\partial\Omega$.
The data $\bs f$ and $\bs g$ are chosen such that the velocity $\bs u$ and the pressure $p$ are defined
as $\bs u(x,y) = (\sin(\pi x)\cos(\pi y), -\sin(\pi y)\cos(\pi x))^T$ and $p(x,y) = \sin(\pi x)\cos(\pi y)$ that defines an enclosed flow. Velocity $\bs u$ and pressure $p$ are discretized on two different
sets of Halton nodes $\Omega_{h,\bs u}\subset\bar\Omega$, and $\Omega_{h,p}\subset \Omega_{h,\bs u}\cap
\Omega$ chosen such that the
stability condition established in Westermann \etal~\cite{westermann2025stability} is satisfied. 
In order to apply Schwarz's method we need to decompose both $\Omega_{h,\bs u}$ and $\Omega_{h,p}$.
In order to achieve stability we ensure that each subdomain is finer discretized for the velocity than for the pressure,
and discrete subsets of nodes for both velocity and pressure
exhibit a sufficient overlap. As before, we determine experimentally the minimum size of the subdomains that ensures 
convergence of the Schwarz iteration. The corresponding minimum number $N_{min}^{\bs u}$ of velocity nodes in the subdomains
is reported below as reference, whereas we keep the quotient between the number of pressure and velocity nodes at
approximately $1/4$  as suggested in~\cite{westermann2025stability}.

\begin{figure}[p]
	\hspace{24pt}	
	\begin{minipage}{0.49\textwidth}
		\includegraphics[width = \figureScalee]{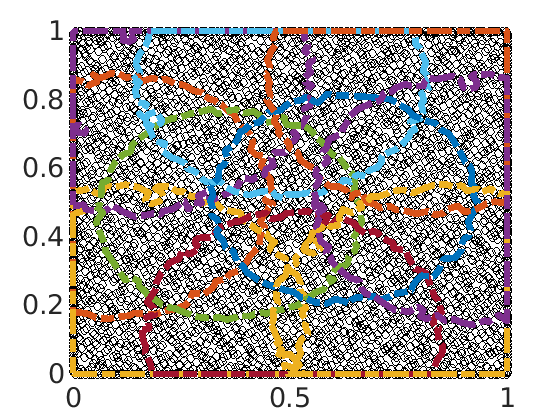}
	\end{minipage}
	\hfill
	\begin{minipage}{0.49\textwidth}
		\includegraphics[width = \figureScalee]{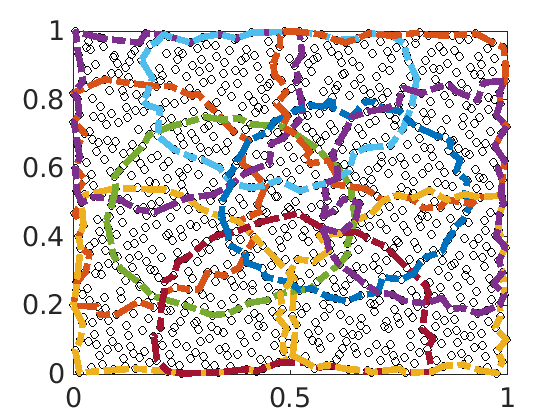}
	\end{minipage}
	\caption{Decomposition example for the Stokes 2D for velocity (left) and pressure (right) with $N_{DEC} = 10$ and 
		$N^{\bs u} =2N_{min}^{\bs u} = 1118$.}
	\label{fig:Example_2Nmin_Stokes}
\end{figure}

\begin{figure}[p]
	\centering
	\begin{subfigure}{0.49\textwidth}
		\centering
		\includegraphics[width = \figureScale]
		{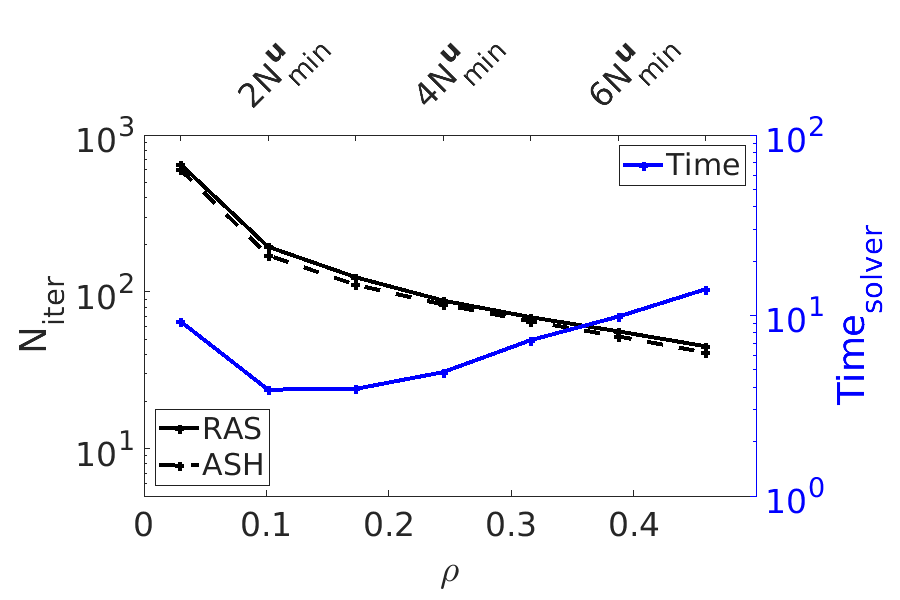}
		\caption{discontinuous PoU}
	\end{subfigure}
	\hfill
	\begin{subfigure}{0.49\textwidth}
		\centering
		\includegraphics[width = \figureScale]
		{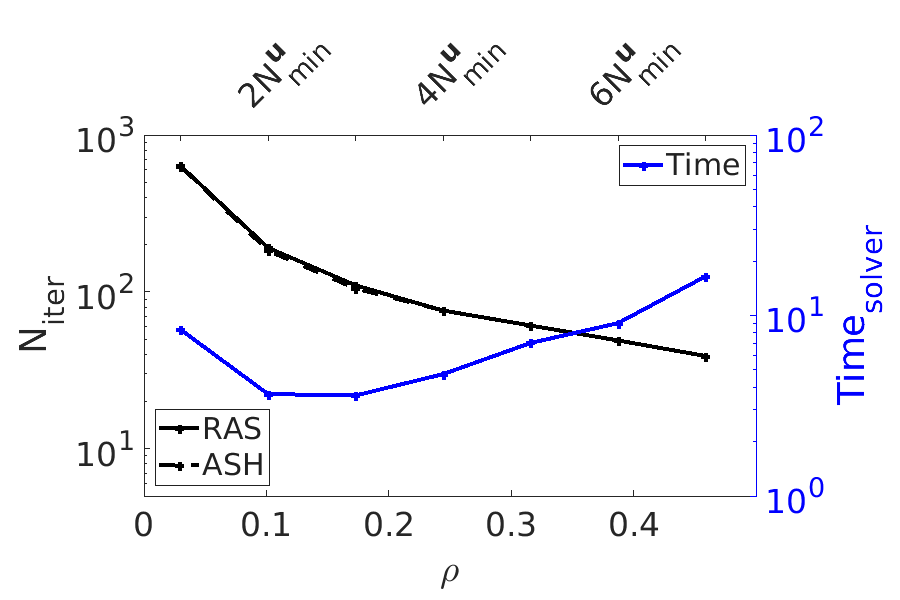}
		\caption{continuous PoU}
	\end{subfigure}
	\begin{subfigure}{0.49\textwidth}
		\centering
		\includegraphics[width = \figureScale]
		{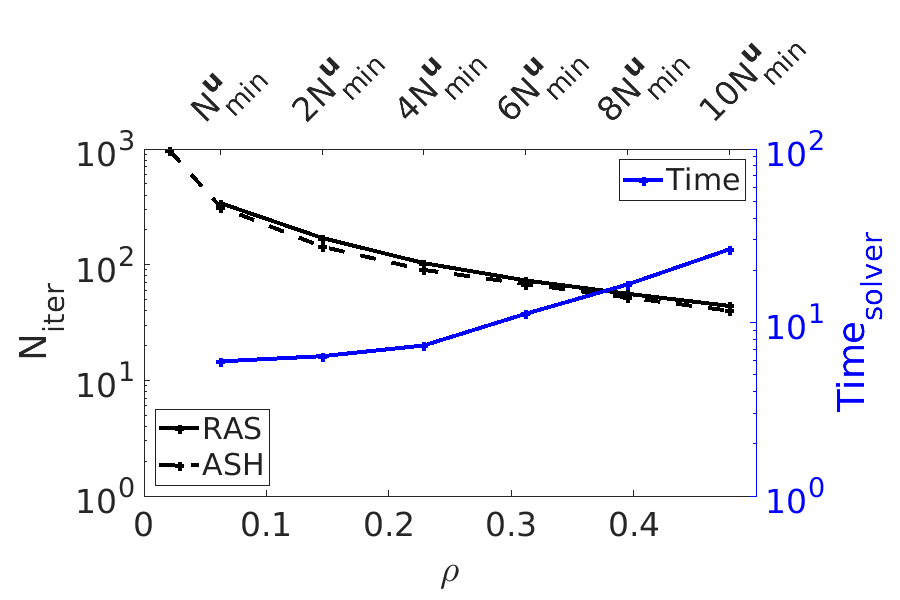}
		\caption{discontinuous PoU}
	\end{subfigure}
	\hfill
	\begin{subfigure}{0.49\textwidth}
		\centering
		\includegraphics[width = \figureScale]
		{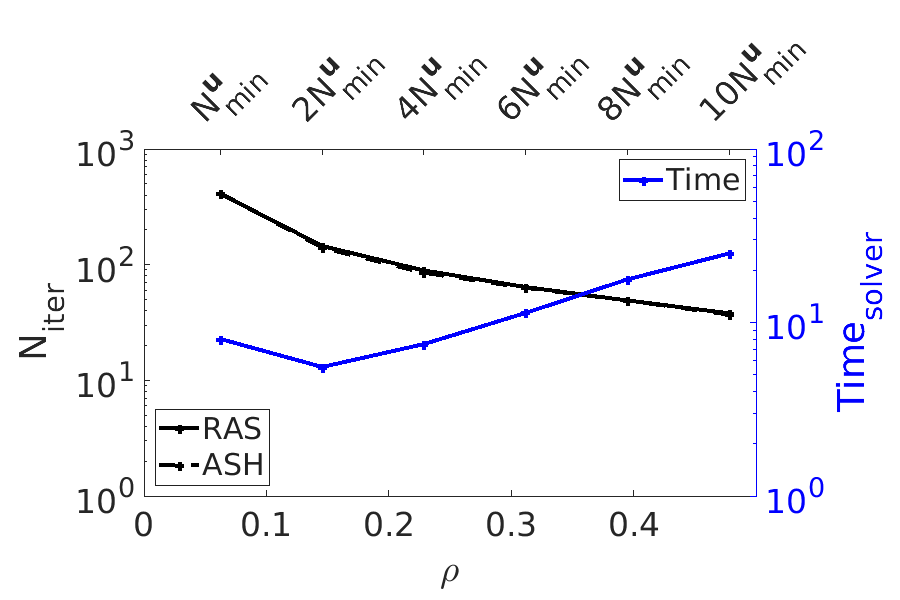}
		\caption{continuous PoU}
	\end{subfigure}
	\begin{subfigure}{0.49\textwidth}
		\centering
		\includegraphics[width = \figureScale]
		{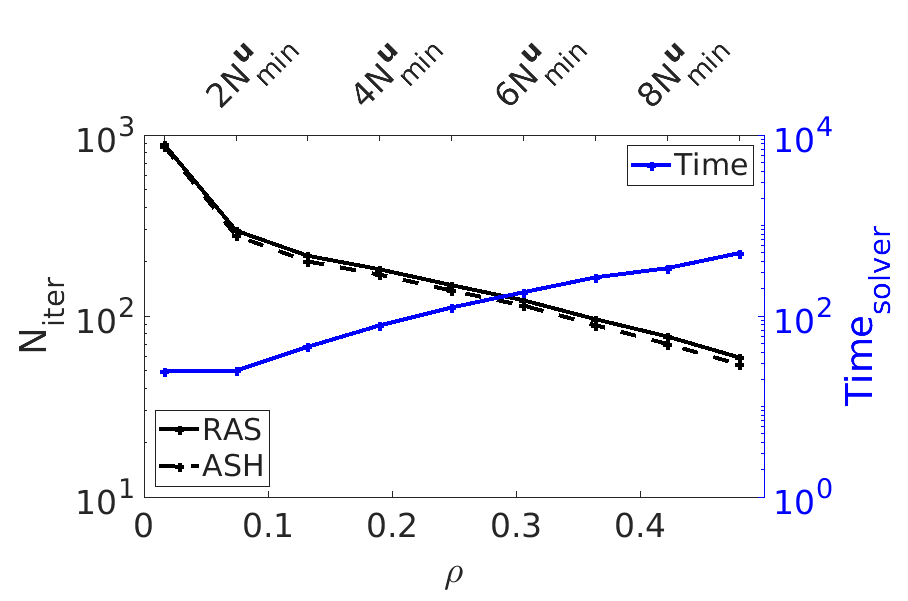}
		\caption{discontinuous PoU}
	\end{subfigure}
	\hfill
	\begin{subfigure}{0.49\textwidth}
		\centering
		\includegraphics[width = \figureScale]
		{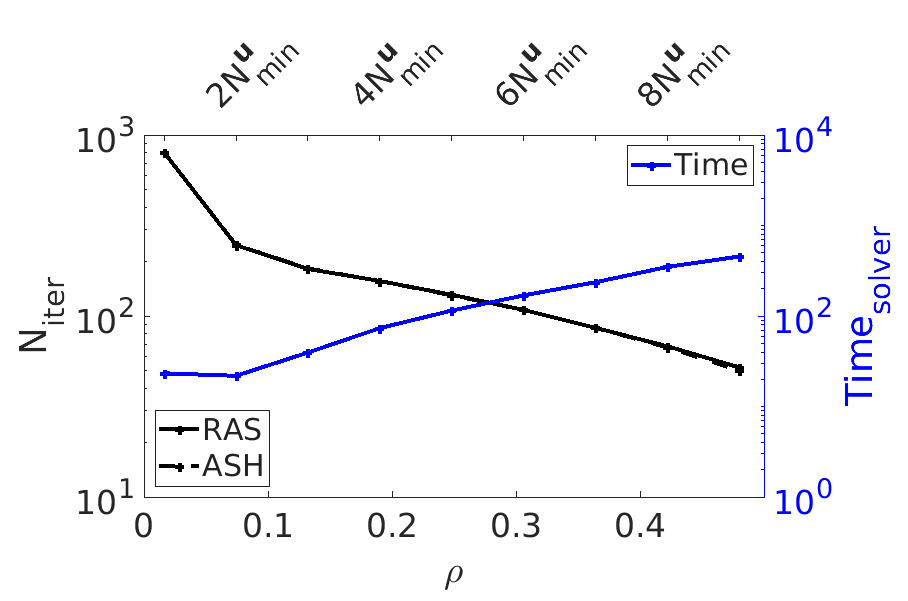}
		\caption{continuous PoU}
	\end{subfigure}
	\caption{Stokes 2D: Comparison of the number of iterations.
	(a)--(b): $N_{max}^{\bs u}=4225$, $N_{DEC} = 25$, $N_{min}^{\bs u}=290$, $N_{dis}^{\bs u}=169$; 
	(c)--(d): $N_{max}^{\bs u}=4225$, $N_{DEC} = 50$, $N_{min}^{\bs u}=172$, $N_{dis}^{\bs u}=85$; 
	(e)--(f): $N_{max}^{\bs u}=16641$, $N_{DEC} = 25$, $N_{min}^{\bs u}=926$, $N_{dis}^{\bs u}=666$.
	}
	\label{fig:Stokes2D_Test}
\end{figure}

Fig.~\ref{fig:Example_2Nmin_Stokes} exemplifies the decomposition obtained for $N^{\bs u}=2N_{min}^{\bs u}$ velocity nodes in
the subdomains, with a near-optimal overlap in our tests, while Fig.~\ref{fig:Stokes2D_Test} provides a  comparison of
the number of iterations and timing for RAS and ASH with dPoU and cPoU. The plots are only shown for $\rho\in[0,0.5]$ as
larger overlaps are of no interest due to a high computational cost. 

We observe that ASH is preferred over RAS, however
with only a small difference. The timing grows faster than in the case of Poisson 2D after passing the minimum 
at $2N_{min}^{\bs u}$, which correlates with the higher cost of the solution by the direct method, 
with the rate $\mathcal{O}(N^{2.8})$ for an LU-decomposition, and $\mathcal{O}(N^{2})$
for solving the system at each iteration. This is due to the complex saddle-point structure of the Stokes equations, which
has a higher bandwidth and is more susceptible to fill-in effects. 

%%%%%%%%%%%%%%%%%%%%
%%%%% Conclusion %%%
%%%%%%%%%%%%%%%%%%%%
\section{Conclusion}\label{Conclusion}
In this paper, we studied the numerical behavior of the additive algebraic Schwarz method  
applied to boundary value problems discretized by the meshless RBF-FD method on irregular nodes. 
We have seen that the Schwarz method can be executed for elliptic and saddle point problems 
on overlapping non-structured decompositions complying with 
the simplicity of the meshless approach.

Our numerical experiments suggest that disjoint subdomains do not work well in this setting, 
and for 2D a near-optimal cost in our experiments is achieved when the overlap is such that 
the number of nodes in the subdomains is about twice 
the minimum number that leads to convergence, meaning a substantial but not too high an overlap, as
illustrated in Fig.~\ref{fig:Example_Nmin} and \ref{fig:Example_2Nmin_Stokes}.

If a discontinuous PoU is used, then the ASH version of the additive algebraic Schwarz method is preferred over RAS 
since the method converges in fewer iterations, while the cost per iteration is comparable. For a
continuous PoU there is a slight advantage of ASH over RAS in 2D examples.

These observations suggest that ASH-type Schwarz method with a small but substantial overlap is a good candidate for a 
smoother in the geometric meshless multigrid methods. 

\vspace{-12pt}

\bibliographystyle{spmpsci}
\bibliography{literature}

@article{cai1999restricted,
	title={A restricted additive {S}chwarz preconditioner for general sparse linear systems},
	author={Cai, Xiao-Chuan and Sarkis, Marcus},
	journal={Siam Journal on Scientific Computing},
	volume={21},
	number={2},
	pages={792--797},
	year={1999},
	publisher={SIAM}
}

@incollection{St-Cyr2007optimized,
	title={Optimized restricted additive {S}chwarz methods},
	author={St-Cyr, Amik and Gander, Martin J and Thomas, Stephen J},
	booktitle={Domain decomposition methods in science and engineering XVI},
	pages={213--220},
	year={2007},
	publisher={Springer}
}

@article{efstathiou2003restricted,
	title={Why Restricted Additive {S}chwarz Converges Faster than Additive {S}chwarz},
	author={Efstathiou, Evridiki and Gander, Martin J},
	journal={BIT Numerical Mathematics},
	volume={43},
	number={5},
	pages={945--959},
	year={2003},
	publisher={Springer}
}

@book{fornberg2015primer,
	title={A primer on radial basis functions with applications to the geosciences},
	author={Fornberg, Bengt and Flyer, Natasha},
	year={2015},
	publisher={SIAM}
}

@article{gander2008schwarz,
	title={Schwarz methods over the course of time},
	author={Gander, Martin J and others},
	journal={Electron. Trans. Numer. Anal},
	volume={31},
	number={5},
	pages={228--255},
	year={2008}
}

@inproceedings{gander2018does,
	title={Does the Partition of Unity Influence the Convergence of {S}chwarz Methods?},
	author={Gander, Martin J},
	booktitle={International Conference on Domain Decomposition Methods},
	pages={3--15},
	year={2018},
	organization={Springer}
}

@article{john2001higher,
	title={Higher-order finite element discretizations in a benchmark problem for incompressible flows},
	author={John, Volker and Matthies, Gunar},
	journal={International Journal for Numerical Methods in Fluids},
	volume={37},
	number={8},
	pages={885--903},
	year={2001},
	publisher={Wiley Online Library}
}

@article{schwarz1870original,
	title={Ueber einen {G}renz\"ubergang durch alternirendes {V}erfahren},
	journal ={Vierteljahrsschrift der Naturforschenden Gesellschaft in Zurich},
	author={Schwarz, HA},
	year={1870}
}

@book{toselli2004domain,
	title={Domain Decomposition Methods - Algorithms and Theory},
	author={Toselli, Andrea and Widlund, Olof},
	volume={34},
	year={2004},
	publisher={Springer Science \& Business Media}
}

@article{vanka1986block,
	title={Block-Implicit Multigrid Solution of {N}avier-{S}tokes Equations in Primitive Variables},
	author={Vanka, S Pratap},
	journal={Journal of Computational Physics},
	volume={65},
	number={1},
	pages={138--158},
	year={1986},
	publisher={Elsevier}
}

@article{westermann2025stability,
	title={Stability and Accuracy of a Meshless Finite Difference Method for the {S}tokes Equations},
	author={Westermann, Alexander and Davydov, Oleg and Sokolov, Andriy and Turek, Stefan},
	journal={International Journal for Numerical Methods in Engineering},
	volume={126},
	number={3},
	year={2025},
	publisher={Wiley Online Library}
}

@misc{mFDlab,
  author={Oleg Davydov}, 
  title={{mFDlab}: A laboratory for meshless finite difference ({mFD}) methods, \url{https://bitbucket.org/meshlessFD/mfdlab}}, 
  year={2020}, 
}

@article{wobker2009numerical,
  title={Numerical studies of {V}anka-type smoothers in computational solid mechanics},
  author={Wobker, Hilmar and Turek, Stefan},
  journal={Advances in Applied Mathematics and Mechanics},
  volume={1},
  number={1},
  pages={29--55},
  year={2009}
}

\end{document}